\documentclass[]{spie}  

\usepackage{amsmath,amsfonts,amssymb}
\usepackage{graphicx}
\usepackage[colorlinks=true, allcolors=blue]{hyperref}
\usepackage{empheq}

\newcommand{\cLtwo}{{\mathcal L}_2 (0,L)}

\usepackage{pifont}
\usepackage[tableposition=top]{caption}
\def\mX{{\mathbb X}}

\newtheorem{rmk}{Remark}[section]
\newtheorem{thm}{Theorem}[section]

\newtheorem{prop}{Proposition}[section]

\newtheorem{lem}{Lemma}[section]

\newcommand{\mc}{\mathcal}
\newcommand{\mb}{\mathbf}

\def\e1{{\varepsilon_{1}}}
\def\b1{{\beta_{11}}}
\def\bp3{{\beta_{33}}}
\def\ep3{{\varepsilon_{3}}}

\MHInternalSyntaxOn
\providecommand*\phantomword[3][c]{%
\mathchoice
{\MT_phantom_word:NNnn #1\displaystyle {#2}{#3}}%
{\MT_phantom_word:NNnn #1\textstyle {#2}{#3}}%
{\MT_phantom_word:NNnn #1\scriptstyle {#2}{#3}}%
{\MT_phantom_word:NNnn #1\scriptscriptstyle {#2}{#3}}%
}
\def\MT_phantom_word:NNnn #1#2#3#4{%
\@begin@tempboxa\hbox{$\m@th#2#4$}%
\setlength\@tempdima{\widthof{$\m@th#2#3$}}%
\hbox{\hb@xt@\@tempdima{\csname bm@#1\endcsname}}%
\@end@tempboxa}
\MHInternalSyntaxOff

\title{Modeling and stabilization results for a charge or current-actuated active constrained layer (ACL) beam model with the electrostatic assumption \footnote{Article submitted to: Smart Structures NDE 2016: Active and Passive Smart Structures and Integrated Systems X, Piezo-based Materials and Systems II, Proc. of SPIE.}}

\author[a]{Ahmet \"Ozkan \"Ozer}
\affil[a]{Department of Mathematics \& Statistics, University of Nevada, 1661 N. Virginia St., Reno, NV 89503, USA}

\authorinfo{Further author information: (Send correspondence to Ahmet \"Ozkan \"Ozer)\\ E-mail: aozer@unr.edu, Telephone: 1-775-7846774
}

\begin{document}
\maketitle

\begin{abstract}
An infinite dimensional model for a three-layer active constrained layer (ACL) beam model, consisting of a piezoelectric elastic layer at the top and an elastic host layer at the bottom constraining a viscoelastic layer in the middle, is obtained  for clamped-free boundary conditions by using a thorough variational approach. The Rao-Nakra thin compliant layer approximation is adopted to model the sandwich structure, and the electrostatic approach (magnetic effects are ignored) is assumed for the piezoelectric layer. Instead of the voltage actuation of the piezoelectric layer, the piezoelectric layer is proposed to be activated by a charge (or current) source. We show that, the closed-loop system with all mechanical feedback is shown to be uniformly exponentially stable. Our result is the outcome of the  compact perturbation  argument and a unique continuation result for the spectral problem which relies on the multipliers method. Finally, the modeling methodology of the paper is generalized  to the multilayer ACL beams, and the uniform exponential stabilizability result is established analogously.
\end{abstract}


\section{Introduction}
\begin{figure}
\begin{center}
\includegraphics[height=5.5cm]{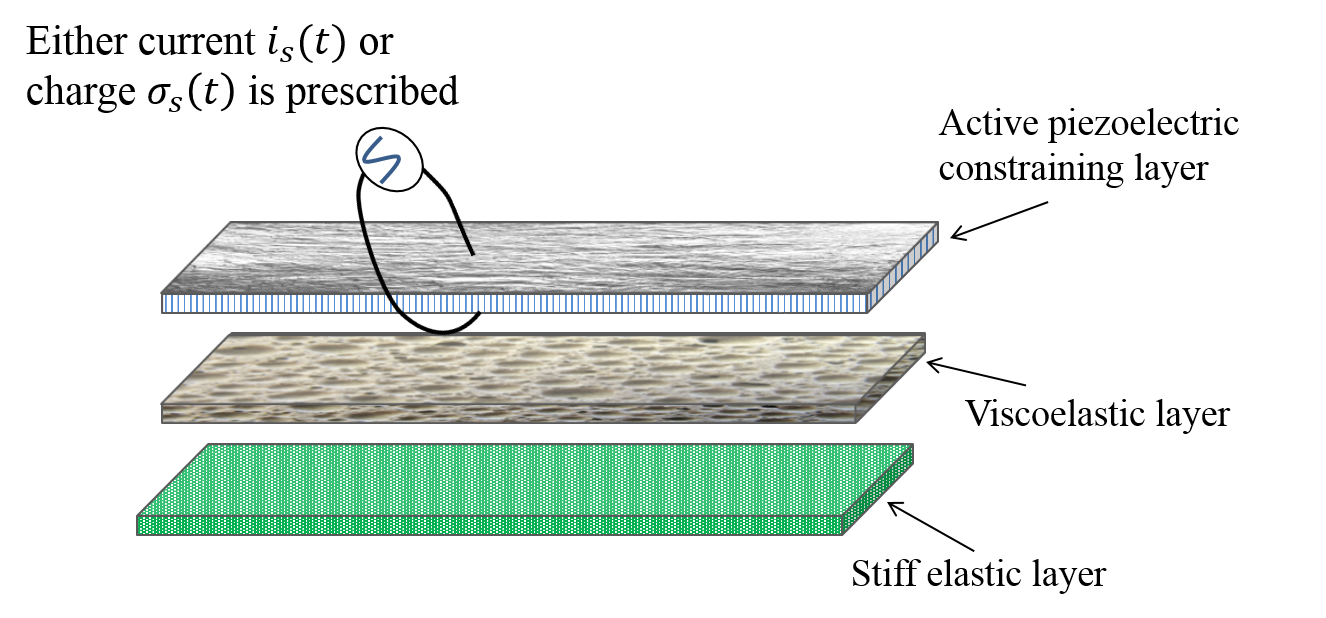}    
\caption{For a current or charge-actuated ACL beam, when   $i_s(t)$ or $\sigma_s(t)$ is supplied to the electrodes of the piezoelectric layer, an electric field is created between the electrodes, and therefore the piezoelectric beam either shrinks or extends, and this causes the whole composite stretch and bend.}  
\label{ACL}                                 
\end{center}                                 
\end{figure}

A three-layer actively constrained layer (ACL) beams is an elastic beam consisting of a stiff elastic beam, a piezoelectric beam, and a viscoelastic beam which creates passive damping in the structure. The piezoelectric beam itself an elastic beam with electrodes at its top and bottom surfaces, insulated at the edges (to prevent fringing effects), and connected to an external electric circuit (see Fig. \ref{ACL}).  Piezoelectric structures are widely used in  in civil, aeronautic and space space structures due to their small size and high power density. They convert mechanical energy to electric energy, and vice versa. ACL composites involve a piezoelectric layer and therefore utilize the benefits of it. Modeling these composites requires better understanding of the modeling of  piezoelectric layer since the they  are generally actuated through the piezoelectric layer. 

There are mainly  three ways to electrically actuate piezoelectric materials: voltage, current or charge. Piezoelectric materials have been traditionally activated by a voltage source \cite{Baz,Cao-Chen,Rogacheva,Smith,Stanway}, and the references therein. However, it has been observed that the type of actuation changes the controllability characteristics of the host structure \cite{O-M3}.
Charge or current controlled piezoelectric beams show \%85 less hysteresis (electrical nonlinearity) than the voltage actuated ones, i.e. see \cite{Cao-Chen,Review,Main1,Main2,M-F}, and the references therein. In the case of voltage and charge actuation, the underlying control operator is unbounded in the energy space whereas in the case of current actuation (including magnetic effects) the control operator is bounded \cite{O-M3}.

Accurately modeling an  ACL beam  also requires understanding of how  sandwich structures are modeled and how the interaction of the elastic layers are established. A three-layer sandwich beam consists of stiff outer layers and a  viscoelastic core layer. The core layer is supposed to deform only in transverse shear. The bending is uniform for the whole composite. Many sandwich beam models have been proposed in the literature, i.e., see \cite{Ditaranto,Mead,Rao,Sun,Trindade,Yan}, and the references therein.  These models mostly differ by the  assumptions on the viscoelastic layer. For example, the Mead-Marcus type models  disregard the effects of longitudinal and rotational inertias \cite{Mead}, and the Rao-Nakra type models preserve these effects since it is  noted that these inertia effects are expected to have considerable importance especially at the high frequency modes for sandwich beams \cite{Rao}.

      The models for the ACL beams proposed in the literature mostly use the  sandwich beam assumptions for the interactions of the layers, i.e. \cite{Stanway,Trindade}, and references therein.  The massive majority of these models are actuated by a voltage source \cite{Baz,F}. Moreover, the longitudinal vibrations were not taken into account. Only the bending of the composite is studied.

      In this paper, we use the Rao-Nakra thin-complaint layer approximation \cite{Hansen3} for which the longitudinal and rotational inertia terms are kept. We obtain two different models, one for the charge actuated ACL beam, another one for the current actuated ACL beam. We assume the electrostatic approach together with the quadratic-through thickness electric potential so that  the induced potential effect is taken into account.
 We show that the proposed model  with charge actuation can be written in semigroup formulation (state-space formulation), and it is well-posed in the natural energy space.  The biggest advantage of our models is being able to study controllability and stabilizability problems for ACL beams in the infinite dimensional setting. We show that the charge actuated model can be uniformly exponentially stabilized by choosing appropriate feedback controllers  which are all mechanical and collocated in this paper due to the electrostatic assumption. This type of feedback controllers avoids the well-known \emph{spillover} effect. Our stability result relies on the compact perturbation  argument and a unique continuation result for the spectral problem by the multipliers technique. A similar argument was used earlier in \cite{O-Hansen3,O-Hansen4} for different boundary conditions.

\section{Modeling Active-Constrained Layer
 (ACL) beams} \label{modeling}
\noindent \textbf{I- Charge or current-actuation:}

The ACL beam is a composite consisting of three layers that occupy the
region $\Omega=\Omega_{xy}\times (0, h)=[0,L]\times [-b,b] \times (0,h)$ at equilibrium where $\Omega_{xy}$ is a smooth bounded domain in the plane.
The total thickness $h$ is assumed to be small in comparison to the dimensions of $\Omega_{xy}$.
The beam  consists of a stiff layer, a compliant layer, and a piezoelectric layer, see Figure \ref{ACL}.

Let $0=z_0<z_1<z_2<z_3=h, $ with
 $$h_i=z_i-z_{i-1}, \quad i=1,2,3.$$
We use the rectangular coordinates $(x,y)$ to denote points in $\Omega_{xy},$ and  $(X, z)$ to denote points in $\Omega = \Omega^{\rm s} \cup \Omega^{\rm ve} \cup \Omega^{\rm p} $, where $\Omega^{\rm s}, \Omega^{\rm ve},$ and $\Omega^{\rm p}$ are the reference configurations of the stiff, viscoelastic, and piezoelectric layers, respectively, and they are defined by
\begin{eqnarray}
\nonumber &&\Omega^{\rm s}=\Omega_{xy}\times (z_0,z_1),\quad  \Omega^{\rm ve}=\Omega_{xy}\times (z_1,z_2), \quad  \Omega^{\rm p}=\Omega_{xy}\times (z_2,z_3).
\end{eqnarray}

 Define
$$\vec \psi=[\psi^1,\psi^2, \psi^3]^{\rm T}, \quad \vec \phi=[\phi^1,\phi^2,\phi^3]^{\rm T}, \quad \vec v=[v^1, v^2, v^3]^{\rm T}$$ where
\begin{eqnarray}
\label{defs1} &&  \psi^i=\frac{u^i-u^{i-1}}{h_i}, \quad \phi^i= \psi^i + w_x, \quad  v^i= \frac{u^{i-1}+u^i}{2}, \quad i = 1, 2, 3,\\
\label{defs3}&&\phi^1=\phi^3=0, \quad \quad \psi^1=\psi^3=-w_x, \quad \phi^2=\psi^2+ w_x
\end{eqnarray}

Let $T_{ij}$ and $ S_{ij}$   denote the stress and strain tensors for $i, j = 1, 2, 3$, respectively.
The constitutive equations for the piezoelectric layers are

\begin{eqnarray}
 \label{cons-eq50}
   T_{11}^{(3)}=\alpha^3 S^{(3)}_{11}-\gamma D_3, ~~   D_1=\varepsilon_{11} E_1, ~~  D_3=\gamma  S_{11}^{(3)}+\varepsilon_{33} E_3,
\end{eqnarray}
where $D, E, \alpha, \gamma,$  and $\varepsilon$ are electrical displacement vector, electric field intensity vector, elastic stiffness coefficient, piezoelectric coefficient,   permittivity coefficient, impermittivity coefficient, and and for the middle and the stiff layers are  given as the following:
\begin{eqnarray}
\label{cons-eq60}
 &  T_{11}^{(i)}=\alpha^i S_{11}^{(i)},\quad T_{13}^{(i)}= 2G_{2} S_{13}^{(i)}, \quad i=1,2&
\end{eqnarray}
where $G_2$ is the shear modulus of the second layer. Since we don't allow shear in the stiff layer we indeed have $T_{13}^{(i)}=0,$ $i=1,3.$ The strain components for the middle layer are
\begin{eqnarray}  \label{strains1} && S_{11}^{(2)}=\frac{\partial v^2}{\partial x}- (z-\hat z_i) \frac{\partial \psi^2}{\partial x}, \quad ~ S_{13}^{(2)}=\frac{1}{2}\left(\psi^2+ w_x\right)=\frac{1}{2}\phi^2,~
\end{eqnarray}
and for the piezoelectric and stiff layers are given by
\begin{eqnarray}
 \label{strains2} & S_{11}^{(i)}=\frac{\partial v^i}{\partial x}- (z-\hat z_i) \frac{\partial^2 w}{\partial x^2},\quad S^{(i)}_{13}=0,\quad i=1,3.&
\end{eqnarray}
For further information about elastic and piezoelectric constants refer to \cite{O-M1}.

\textbf{Electrostatic modeling:} In contrast to the dynamic modeling \cite{O-M3,Ozkan3}, we assume the electrostatic assumption, and therefore the set of Maxwell's equations with the appropriate mechanical boundary conditions at the edges of the beam  (the beam is clamped, hinged, free, etc.) is partially used. In other words, the magnetic field, $\frac{\partial \dot B}{\partial t},$ and $\frac{\partial\dot D}{\partial t}$ are ignored. Depending on the type of actuation the charge density $\sigma_s,$ the current density $i_s,$ or voltage $V$ is prescribed at the electrodes, i.e.  on the faces $\partial\Omega^p.$  In this paper we consider  the charge (and current).  The voltage actuation case is handled in details in \cite{Ozkan2}.

 By the Faraday's law,  there exists a scalar electric potential $\varphi$  such that \begin{equation}
\label{imp1} \quad E=-\nabla\varphi.
\end{equation}
   Henceforth, to simplify the notation, $x=x_1$ and $z=x_3.$

 The linear through-thickness assumption of the electric potential $\varphi(x,z)=\varphi^0(x)+z\varphi^1(x),$ which is a common assumption in many papers,  completely ignores the induced potential effect since $\varphi$ is completely known as a function of voltage \cite{Tabesh}. Therefore  we use a quadratic-through thickness potential distribution to take care this effect:
 \begin{eqnarray}
 \label{scalarpot} && \varphi(x,z)=\varphi^0(x)+z \varphi^1(x)+ \frac{z^2}{2} \varphi^2(x)
 \end{eqnarray}
By (\ref{imp1})
\begin{eqnarray}
\nonumber &&E_1=   -\left(({\varphi}^0)_x+z (\varphi^1)_x+ \frac{z^2}{2} (\varphi^2)_x\right),\\
\label{e3} &&E_3=-\left( \varphi^1 + z \varphi^2 \right).
\end{eqnarray}

\noindent {{\bf Work done by the external forces:}} The work done by the electrical external forces (as in \cite{Lee,O-M3}) is
\begin{eqnarray}
\nonumber &&\mb{W}^{e}=  \int_{\partial \Omega^{\rm p}} \left(-\tilde\sigma_s~\phi + \tilde i_s\cdot  A\right)   ~ dX\\
 \nonumber &=& \int_0^L  \left(-\tilde \sigma_s \left(\phi(h_3/2)-\phi(-h_3/2)\right) + \tilde i_s^1\left(A(h_3/2)-A(-h_3/2)\right)\right)  ~dx
\end{eqnarray}
 where  $A$ is the magnetic vector potential, $i_s(x,t)$ is the current density, and $\sigma_s(t)$ is the charge density at the electrodes. Depending on the desired type of actuation, either $i_s(t)$ or $\sigma_s(t)$ is chosen to be zero. In the case of fully dynamic and quasi-static approach, $A$ is assumed to be nonzero. However, the electrostatic approach assumes that $A\equiv 0,$ and therefore the current density $i_s(t)$ does not show up in $\mb W_e$ through the variational approach. From this point on, we continue only with $\sigma_s(t)$ yet $i_s(t)$ may come into the play later through a circuit equation. In fact, one can derive the current controlled ACL beam through the voltage or charge controlled models by adding an additional circuit equation, see further the models Section \ref{models}. Note that this does not happen so in the case of quasi-static or fully dynamic cases. There has to be a compatibility condition to be satisfied \cite{O-M3}. Note also that the term  $\int_{\partial\Omega^p} \dot D_i n_i dS $ is the current flowing out of the electrode.  Thus,
 \begin{eqnarray}
\label{work}  &&\mb{W}^{e}=  \int_{\partial {\Omega^{\rm p}}} -\tilde\sigma_s~\phi    ~ dX \quad {{\rm (charge ~actuation)}}
\end{eqnarray}

Assume that the beam is subject to a  distribution of forces $(\tilde g^1, \tilde g^3, \tilde g)$ along its edge $x=L.$ In parallel to \cite{Hansen3}, then the total work done by the mechanical external forces is
\begin{eqnarray*}
\label{work-done}  \mb{W}^m&=&   g^1 v^1 (L)  + g^3 v^3 (L) +  g w(L)-M w_x(L).
\end{eqnarray*}
where $M$ is the total applied moment at $x=L.$

Now we use the constitutive equations (\ref{cons-eq50})- (\ref{strains2}), and (\ref{scalarpot})-(\ref{e3}) to find 
\begin{eqnarray}
\nonumber &&\mb P=\mb{P}^1 + \mb{P}^2=\frac{1}{2}\int_{\Omega^{\rm s} \cup \Omega^{\rm ve}} \sum_{i=1,2}\left(T_{11}^{(i)}S_{11}^{(i)}+ T_{13}^{(i)}S_{13}^{(i)}\right)~dX \\
\nonumber &&= \frac{1}{2} \int_0^L \left[{\alpha}^1 h_1\left((v^1_x)^2+ \frac{h_1^2}{12}w_{xx}^2\right)+{\alpha}^2 h_2\left((v^2_x)^2+ \frac{h_2^2}{12}(\psi_x)^2\right)+ G_2h_2 (\phi^2)^2\right]~dx,
\end{eqnarray}
\begin{eqnarray}
\nonumber \mb{E}^3-\mb{P}^3&=&\frac{1}{2}\int_{\Omega^{\rm p}} \left(D_1 E_1+D_3 E_3- T_{11}^{(3)}S_{11}^{(3)}-T_{13}^{(3)}S_{13}^{(3)}\right)~dX \\
\nonumber &=& \frac{1}{2} \int_0^L \left[-{\alpha}^3 h_3\left((v^3_x)^2+ \frac{h_3^2}{12}(w_{xx})^2\right) - 2{\gamma}h_3  \left( \varphi^1 v^3_x-\frac{h_3^2}{12}w_{xx} \varphi^2\right)\right.\\
 \label{E-P} && + \left.{\e1}h_3 \left( (\varphi^0_x)^2+\frac{h_3^2}{12} (\varphi^1_x)^2+\frac{h_3^2}{24} (\varphi^0)_x( \varphi^2)_x + \frac{h_3^4}{320} (\varphi^2_x)^2   \right)  + {\ep3}h_3 \left( (\varphi^1)^2+\frac{h_3^2}{12} (\varphi^2)^2   \right)  \right)~dx
\end{eqnarray}
and $\mb K=\mb K^1+\mb K^2+\mb K^3 $ where
 \begin{eqnarray}
\label{KK1} \mb{K}^i&=&  \frac{\rho_i h_i}{2} \int_0^L \left((\dot v^i)^2 + \dot w^2 +\frac{h_i^2}{12}\dot w_x^2\right)~dx,\quad i=1,3,\\
\label{KK2} \mb{K}^2&=&  \frac{\rho_2 h_2}{2} \int_0^L \left((\dot v^2)^2 + (\dot \psi^2)^2+ \dot w^2 \right)~dx.
\end{eqnarray}


\subsection{Variational Principle \& Equations of Motion}
By using (\ref{defs1})-(\ref{defs3}), the variables $\{v^2, \phi^2, \psi^2\}$ can be written as the functions of $\{v^1, v^3, w\}.$
 For that reason, we choose  $\{v^1, v^3, w,   \phi^1\}$ as the state variables.

To model charge or current-actuated ACL beams with magnetic effects, we use the following Lagrangian \cite{Lee,O-M3}
\begin{eqnarray}\label{Lag}  \mb{ L}= \int_0^T \left[\mb{K}-(\mb{P}-\mb{E})+ \mb{W}\right]dt\end{eqnarray}
where $\mb P-\mb E=\mb P^1 + \mb P^2 +\mb P^3 -\mb E^3 $ is called electrical enthalpy where $\mb W=\mb W^m+\mb W^e$ is the total work done by the mechanical and electrical external forces, respectively. Note that in modeling piezoelectric beams by voltage-actuated electrodes we use a modified Lagrangian  \cite{O-M1}. Let $H=\frac{h_1 + 2h_2+h_3}{2}.$ We assume that the ACL beam is clamped at $x=0$ and free at $x=L$.

The application of Hamilton's principle, setting the variation of Lagrangian $\mb L$ in  (\ref{Lag}) with respect to the all kinematically  admissible displacements of the state variables $\{v^1, v^3, w,   \phi^1\}$ to zero, yields a strongly coupled equations for the longitudinal and transverse dynamics. It is not easy to study the controllability/stabilizability  properties. For this reason we study the approximated model.

\subsection{Rao-Nakra model assumptions}
\label{models}
In this section we derive the  approximated Rao-Nakra type ACL beam model assumptions by assuming that the compliant layer is thin, i.e. $\rho_2, \alpha^2\to 0.$  This variation of the initial model corresponding to thin compliant layer is  described and shown to be a regular perturbation of the initial model \cite{Hansen3}.   As well, this
approximation retains the potential energy of shear and transverse kinetic energy. We obtain the following model

\begin{eqnarray}
\label{reduced-nog} \left\{
  \begin{array}{ll}
   m \ddot w - K_1 \ddot{w}_{xx}+K_2w_{xxxx} -  H G_2   \phi^2_x=0&\\
  \rho_1h_1 \ddot v^1  -\alpha^1 h_1 v^1_{xx}  - G_2  \phi^2 = 0,   & \\
  \rho_3h_3 \ddot v^3       -\alpha^3 h_3 v^3_{xx}  + G_2 \phi^2 -{\gamma}h_3 \left((\varphi^1)_x \right) = 0,   & \\
    - \xi(\phi^1)_{xx} +\phi^1 -\frac{{\gamma}}{\ep3} v^3_x= \frac{\sigma_s (t) }{h_3\ep3} & \\
      \phi^2 =\frac{1}{h_2}\left(-v^1+v^3\right) + \frac{H}{h_2}w_x,& ~~~~ (x,t)\in (0,L)\times \mathbb{R}^+
  \end{array} \right.
  \end{eqnarray}
with the natural boundary  conditions at $x=0,L$
\begin{eqnarray}
&&  \begin{array}{ll}
\nonumber   \left\{v^1, v^3, w, w_x \right\}_{x=0}=0, ~~   \left\{\alpha^1 h_1v^1_x\right\}_{x=L}=g^1(t),~~ \left\{ \alpha^3 h_3 v^3_x  +{\gamma}h_3\phi^1\right\}_{x=L}=g^3(t)&\\
  \end{array} \\
      &&  \begin{array}{ll}
\nonumber  \left\{  (\phi^1)_x\right\}_{x=0,L}=0,~~   \left\{K_2 w_{xx}  \right\}_{x=L} =-M(t), ~ \left\{K_1 \ddot w_x -K_2 w_{xxx} + G_2 H \phi^2 \right\}_{x=L} =g(t),& t\in \mathbb{R}^+\\
\end{array} \\
      &&  \begin{array}{ll}
\label{reduced-BC-nog}   (v^1, v^3, w, \dot v^1, \dot v^3, \dot w)(x,0)=(v^1_0,  v^3_0, w_0,  v^1_1, v^3_1,  w_1),& ~~ (x,t)\in (0,L)\times \mathbb{R}^+
  \end{array}
\end{eqnarray}
where $m=\rho_1h_1+ \rho_3 h_3,$ $K_1=\frac{\rho_1 h_1^3}{12} +\frac{\rho_3 h_3^3}{12},$ and $K_2=\frac{\alpha^1 h_1^3}{12}+\frac{\alpha^3 h_3^3}{12}.$

Observe that the last equation in (\ref{reduced-nog}) is an elliptic equation. For that reason first we eliminate $\varphi^1$ in (\ref{reduced-nog}). Defining $D_x^2 \phi = \phi_{xx}$ and its domain
$${\rm Dom} (D_x^2) = \{ \phi \in H^2 (0,L) ,\quad \phi_x (0)=\phi_x (L)=0 \}, $$
and the operator $P_\xi$ by
\begin{eqnarray}\label{Lgamma}P_{\xi}:=\left(-\xi D_x^2+I\right)^{-1}, \quad \xi:=\frac{\e1 h_3^2}{12\ep3}.\end{eqnarray}
It is well-known that ${P_\xi}$ is a non-negative and a compact operator on $\cLtwo$ \cite{O-M3}.  Thus, the system (\ref{reduced-nog})-(\ref{reduced-BC-nog}) is simplified to
\begin{eqnarray}
\label{electro} \left\{
  \begin{array}{ll}
    m \ddot w - K_1 \ddot{w}_{xx}+K_2w_{xxxx} -  H G_2   \phi^2_x=0,&\\
    \rho_1h_1 \ddot v^1  -\alpha^1 h_1 v^1_{xx}  - G_2  \phi^2 =0,   &  \\
  \rho_3h_3 \ddot v^3       -\alpha^3 h_3 v^3_{xx} -\frac{{\gamma^2 h_3}}{{\ep3}}~(P_{\xi} v^3_x)_x  + G_2 \phi^2  = \frac{\gamma\sigma_s(t)}{\ep3 } \delta(x-L),   & \\
     \phi^2 =\frac{1}{h_2}\left(-v^1+v^3\right) + \frac{H}{h_2}w_x,&  ~~ (x,t)\in (0,L)\times \mathbb{R}^+
  \end{array} \right.
  \end{eqnarray}
with the initial and natural boundary  conditions at $x=0,L$ and initial conditions
\begin{eqnarray}
&&  \begin{array}{ll}
\nonumber   \left\{v^1, v^3, w, w_x  \right\}_{x=0}=0, ~  \left\{\alpha^1 h_1v^1_x\right\}_{x=L}=g^1(t),~~  \left\{ \alpha^3 h_3 v^3_x  +\frac{{\gamma^2 h_3}}{{\ep3}}~P_{\xi} v^3_x\right\}_{x=L}=0,&\\
  \end{array} \\
         &&  \begin{array}{ll}
\nonumber  \left\{K_2 w_{xx}  \right\}_{x=L} =-M(t), ~  \left\{K_1 \ddot w_x -K_2 w_{xxx} + G_2 H \phi^2 \right\}_{x=L} =g(t),& t\in \mathbb{R}^+ \\
  \end{array} \\
&&  \begin{array}{ll}
\label{electro-BC}  (v^1,v^3, w, \dot v^1, \dot v^3, \dot w)(x,0) =(v^1_0, v^3_0, w_0, v^1_1, v^3_1, w_1),& x \in [0,L].
  \end{array}
\end{eqnarray}
where $m=\rho_1h_1+\rho_3 h_3$, $H=\frac{h_1+2h_2+h_3}{2h_2},$ and $ \rho_1,\rho_3, \alpha_1, \alpha_3, \gamma, K_1, K_2, G, h_1, h_3$ are....$\sigma(t)$ is the charge density on the electrodes.

\noindent \textbf{II-Voltage actuation:} Instead of current or charge actuation, one can also go with the traditional voltage actuation \cite{Cao-Chen,Rogacheva,Smith}. This case is analyzed in \cite{Ozkan2,Ozkan3}, and the following model is obtained with the same Rao-Nakra assumptions  
\begin{eqnarray}
 \label{d4-non} &&\left\{
  \begin{array}{ll}
   m \ddot w - K_1 \ddot{w}_{xx} +K_2  w_{xxxx} - {G_2 H}  \phi^2_x= 0&\\
   \rho_1h_1\ddot v^1-\alpha^1h_1 v^1_{xx} - G_2  \phi^2 = 0   & \\
    \rho_3h_3\ddot v^3-\alpha_{1}^3 h_3 v^3_{xx} + G_2  \phi^2 = 0   &\\
    \phi^2=\frac{1}{h_2}\left(-v^1+v^3 + H w_x\right),& ~~ (x,t)\in (0,L)\times \mathbb{R}^+
  \end{array} \right.
\end{eqnarray}
with the initial and natural boundary  conditions at $x=0,L$
\begin{eqnarray}
&&  \begin{array}{ll}
\nonumber   \left\{v^1, v^3, w, w_x \right\}_{x=0}=0, ~~   \left\{\alpha^1 h_1v^1_x\right\}_{x=L}=g^1(t),~~\left\{   \alpha^3 h_3 v^3_{x}(L)= -\gamma V(t)\right\}_{x=L}=0&\\
  \end{array} \\
        &&  \begin{array}{ll}
\nonumber   \left\{K_2 w_{xx}  \right\}_{x=L} =-M(t), ~ \left\{K_1 \ddot w_x -K_2 w_{xxx} + G_2 H \phi^2 \right\}_{x=L} =g(t),& t\in \mathbb{R}^+\\
  \end{array}\\
        &&  \begin{array}{ll}
\label{d-son-non}  (v^1, v^3, w, \dot v^1, \dot v^3, \dot w)(x,0)=(v^1_0,  v^3_0, w_0,  v^1_1, v^3_1,  w_1),& ~~ (x,t)\in (0,L)\times \mathbb{R}^+
  \end{array}
\end{eqnarray}
where $V(t)$ is the voltage applied through the electrodes of the piezoelectric layer.
\section{Feedback Stabilization Results and comparison of models}
\label{Section-stab}
In this section,  we compare the results for the stabilization of the models (\ref{reduced-nog})-(\ref{reduced-BC-nog}) and (\ref{d4-non})-(\ref{d-son-non}).
Note that both cases are not much different from each other. The difference arises only from the $P_\xi$ term in (\ref{electro}). This is only because the voltage actuated model ignores the induced effect of the electric field, i.e. see (\ref{scalarpot}). We account for this effect by choosing the electric potential quadratic-through thickness in this paper. This effect, in fact, turns out to make the piezoelectric beam more stiff since $P_\xi$ is a positive and compact  operator. This was first observed in \cite{Hansen}.

Consider the weak solutions of  (\ref{electro})-(\ref{electro-BC}) with the following feedback controllers
\begin{eqnarray}\label{feedback} g^1=-s_1\dot v^1(L,t), ~~\sigma_s(t)=-s_3\dot v^3(L,t),~~M(t)= k_1\dot w_x(L,t),~~ g(t)=k_2 \dot w(L)
\end{eqnarray}
where $s_1,s_3,k_1, k_2>0.$

The well-posedness of the model (\ref{d4-non})-(\ref{d-son-non}) with (\ref{feedback}) is shown in \cite{Ozkan3}. For that reason, we only analyze the well-posedness of (\ref{reduced-nog})-(\ref{reduced-BC-nog}) with (\ref{feedback}).

Define the complex linear spaces
\begin{eqnarray}
\nonumber  \mathrm V&=&\left(H^1_L(0,L)\right)^2 \times H^2_L(0,L),\quad  \mathrm H= \mX^2 \times H^1_L(0,L),\quad \mc{H} = \mathrm V \times \mathrm H.
\end{eqnarray}

The  natural energy associated with  (\ref{d4-non})-(\ref{d-son-non})  is
\begin{eqnarray}
\nonumber  &&\mathrm{E}(t) =\frac{1}{2}\int_0^L \left\{\rho_1 h_1  |\dot v^1|^2 + \rho_3 h_3  |\dot v^3|^2 + m |\dot w|^2 + \alpha^1 h_1 |v^1_x|^2 + \alpha^3 h_3 |v^3_x|^2 \right. \\
\label{Energy-nat-non}  &&\left. + \frac{{\gamma^2 h_3}}{{\ep3}}~(P_{\xi} v^3_x) (\bar v^3)_x + K_1 |\dot w_x|^2 + K_2 |w_{xx}|^2 +G_2 h_2 |\phi^2|^2 \right\}~ dx.
\end{eqnarray}
 This motivates definition of the inner product on $\mc H$
{ \small{\begin{eqnarray}
\nonumber && \left<\left[ \begin{array}{l}
 u_1 \\
 u_2 \\
 u_3\\
 u_4 \\
 u_5 \\
 u_6
 \end{array} \right], \left[ \begin{array}{l}
 v_1 \\
 v_2 \\
 v_3\\
 v_4\\
 v_5\\
 v_6
 \end{array} \right]\right>_{\mc H}= \left<\left[ \begin{array}{l}
 u_4\\
 u_5\\
 u_6
 \end{array} \right], \left[ \begin{array}{l}
 v_4\\
 v_5\\
 v_6
 \end{array} \right]\right>_{\mathrm H}+ \left<\left[ \begin{array}{l}
 u_1 \\
 u_2 \\
 u_3
 \end{array} \right], \left[ \begin{array}{l}
 v_1 \\
 v_2\\
 v_3
 \end{array} \right]\right>_{\mathrm V}\\
\nonumber && =\int_0^L \left\{\rho_1 h_1  u_5 {\dot {\bar v}}_5 + \rho_3 h_3  u_6 {\dot {\bar v}}_6+ \mu h_3  \dot u_7 {\dot {\bar v}_7} + m \dot u_8{\dot {\bar v}_8}+ K_1 (u_8)_x(\bar v_8)_x+ \alpha^1 h_1  (u_1)_{x} (\bar v_1)_x  \right.  \\
\nonumber &&\left. + \alpha^3 h_3  (u_3)_{x} (\bar v_3)_x+ \frac{{\gamma^2 h_3}}{{\ep3}}~(P_{\xi} (u_3)_x) (\bar u_3)_x  + K_2 (u_4)_{xx} (\bar v_4)_{xx}  +\frac{G_2}{h_2}  (-u_1+u_2 + H(u_4)_x)(-\bar v_1+\bar v_2 + H(\bar v_4)_x)\right\}~dx.
 \end{eqnarray}}}
 Since $P_\xi$ is a positive operator and  the term $\|-u_1+u_2 + H(u_4)_x\|_{L^2(0,L)}$ is coercive, see \cite{Hansen3} for the details., $\langle \, , \, \rangle_{\mc H} $ does indeed define an inner product.

  Let $\vec y=(v^1, v^3, w)$ be the smooth solution of the system of (\ref{d4-non})-(\ref{d-son-non}) with (\ref{feedback}).  Multiplying the equations in (\ref{d4-non}) by $\tilde y_1, \tilde y_2 \in H^1_L(0,L)$ and $ \tilde y_3 \in H^2_L(0,L),$ respectively, and integrating by parts yields
\begin{subequations}
\small
 \begin{empheq}[left={\phantomword[r]{0}{}  }]{align}
\nonumber  &\int_0^L \left(\rho_1h_1\ddot v^1 \tilde y_1 + \alpha^1 h_1 v^1_{x} (\tilde y_1)_x - G_2 \phi^2 \tilde y_1\right)dx + s_1\dot v_1(L) y_1(L)= 0,   & \\
  \nonumber &\int_0^L \left(\rho_3h_3\ddot v^3 \tilde y_2 +\alpha^3 h_3 v^3_{x} (\tilde y_2)_x  +  \frac{\gamma^2 h_3}{\ep3} (P_\xi v^3_x) (\tilde y_2)_x+G_2 \phi^2 (\tilde y_2)_x \right)dx + \frac{\gamma s_3}{\ep3} \dot v_3(L) y_2(L) = 0,   & \\
 \nonumber & \int_0^L \left( m \ddot w \tilde y_3 + K_1 \ddot{w}_{x} (\tilde y_3)_x   + K_2 w_{xx} (\tilde y_3)_{xx}- G_2 H \phi^2_x (\tilde y_3)_x \right)dx + k_1 \dot w_x(L) (y_3)_x(L)+  k_2\dot w(L) y_3(L)  =0.&
\end{empheq}
\end{subequations}
Now define the linear operators
\begin{eqnarray}
\nonumber \left<Ay,\psi\right>_{\mathrm V'\times \mathrm V}=(y,\psi)_{\mathrm V \times \mathrm V}, \quad \forall y,\psi\in \mathrm V\\
\nonumber \left<B_0\vec  y, \vec \psi\right>_{\mathrm H'\times \mathrm H}= \left[ \begin{array}{c}
 0_{2\times 1}\\
 k_2y_3(L) \psi_3(L)
 \end{array} \right], \quad \forall \vec y,\vec \psi\in \mathrm H\\
\label{ops}\left<D_0\vec  y, \vec \psi\right>_{\mathrm H'\times \mathrm H}= \left[ \begin{array}{c}
 s_1y_1(L) \psi_1(L) \\
 \frac{\gamma s_3}{\ep3} y_2(L) \psi_2(L) \\
 k_1(y_3)_x(L) (\psi_3)_x(L)
 \end{array} \right], \quad \forall \vec y,\vec \psi\in \mathrm V.
\end{eqnarray}
Let $\mc M : H^1_L (0,L) \to (H^1_L(0,L))'$ be a linear operator defined by
\begin{eqnarray}\label{def-M}\left< \mc M \psi, \tilde \psi \right>_{(H^1_L(0,L))', H^1_L(0,L)}= \int_0^L (m \psi \tilde {\bar \psi} + K_1 \psi_x \tilde {\bar\psi}_x) dx.
\end{eqnarray}
From the Lax-Milgram theorem $\mc M$ and $A$ are  canonical  isomorphisms from $H^1_L(0,L)$ onto  $(H^1_L(0,L))'$ and from $V$ onto $V',$ respectively.
Assume that $A y \in V',$ then we can formulate the variational equation above into the following form
\begin{eqnarray} M\ddot y  + A y +  D_0\dot y +  B_0\dot y=0
\end{eqnarray}
where $M=\left[\rho_1 h_1 I ~~\rho_3 h_3 I~~ \mc M\right]$ is an isomorphism from $\mathrm H$ onto $\mathrm H'.$
Next we introduce the linear unbounded operator by
\begin{equation}
\label{defA}\mc A: {\text{Dom}}(A)\times V\subset \mc H \to \mc H
\end{equation}
where $\mc A= \left[ {\begin{array}{*{20}c}
   O_{3\times 3}  & -I_{3\times 3} \\
       M^{-1}A  &  M^{-1}  D_0  \\
\end{array}} \right]$
with $\nonumber{\rm {Dom}}(\mc A) =  \{(\vec z, \vec {\tilde z}) \in V\times V, A\vec z \in \mathrm V' \}$
and if ${\rm Dom}(\mc A)'$ is the dual of ${\rm Dom}(\mc A)$ pivoted with respect to $\mc H,$  we define the control operators $B$ and $D$
 \begin{eqnarray}
\label{defb_0} \quad B \in \mathcal{L}(\mathbb{C} , {\rm Dom}(\mc A)'), ~ \text{with} ~ B=   \left[ \begin{array}{c} 0_{3\times 1} \\  M^{-1}B_0. \end{array} \right]
 \end{eqnarray}

Writing $\varphi=[v^1, v^3, w, \dot v^1,  \dot v^3,  \dot w]^{\rm T},$ the control  system (\ref{d4-non})-(\ref{d-son-non})  with the feedback controllers (\ref{feedback}) can be put into the  state-space form
\begin{eqnarray}
\label{Semigroup}
\dot \varphi + \mc A \varphi + B\varphi =0, \quad\varphi(x,0) =  \varphi ^0.
\end{eqnarray}

 \begin{lem} \label{skew}The operator $\mc{A}$  defined by (\ref{defA}) is maximal monotone in the energy space $\mc H,$
 and ${\rm Range}(I+\mc A)=\mc H.$
\end{lem}

\textbf{Proof:} Let $\vec z \in {\rm Dom}(\mc A).$ A simple calculation using integration by parts and the boundary conditions yields
{\small{\begin{eqnarray}
\nonumber &&\left<\mc A  \left[ \begin{array}{c}
\vec z_1 \\
\vec z_2 \end{array} \right], \left[ \begin{array}{c}
\vec z_1 \\
\vec z_2 \end{array} \right]\right>_{\mc H\times \mc H} = \left<  \left[ \begin{array}{c}
-\vec z_2 \\
M^{-1} A \vec z_1 \end{array} \right], \left[ \begin{array}{c}
\vec z_1 \\
\vec z_2 \end{array} \right]\right>_{\mc H}= \left<-\vec z_2, \vec z_1\right>_{V\times V}  + \left<M^{-1}\left(A \vec z_1+D_0\vec z_2  \right), \vec z_2\right>_{H\times H}\\
\nonumber &&\quad\quad = -\overline{\left<A \vec z_1, \vec z_2\right>_{V'\times V}}  + \left<A\vec z_1 + D_0z_2, \vec z_2\right>_{H'\times H}.
\end{eqnarray}
}}
Since $\vec z=\left[ \begin{array}{c}
\vec z_1 \\
\vec z_2 \end{array} \right] \in {\rm Dom}(\mc A) $, then $A\vec z _1+ D_0z_2\in V'$ and $\vec z_2 \in V$ so that
\begin{eqnarray*}&&\left<A\vec z_1 +  D_0z_2, \vec z_2\right>_{H'\times H}=\left<A\vec z_1+ D_0z_2, \vec z_2\right>_{V'\times V} =\left<A\vec z_1, \vec z_2\right>_{V'\times V}+\left< D_0z_2, \vec z_2\right>_{V'\times V}.
\end{eqnarray*}
Hence
${\rm Re}\left<\mc A  \vec z, \vec z\right>_{\mc H\times \mc H}=\left< D_0z_2, \vec z_2\right>_{V'\times V}\ge 0.$
We next verify the range condition. Let $\vec z=\left[ \begin{array}{c}
\vec z_1 \\
\vec z_2 \end{array} \right]\in \mc H.$ We prove that there exists a $\vec y=\left[ \begin{array}{c}
\vec y_1 \\
\vec y_2 \end{array} \right]\in \mc {\rm Dom} (\mc A )$ such that $(I+\mc A) \vec y=\vec z.$
A simple computation shows that proving this is equivalent to proving ${\rm Range} (M+A + D_0)=H',$ i.e., for every $\vec f \in H'$   there exists a unique solution $\vec z \in H$ such that $(M+A+D_0)\vec z=\vec f.$
This obviously follows from the Lax Milgram's theorem. $\square$

\begin{prop} The operator $B$ is a monotone compact operator on $\mathrm H.$
\end{prop}

\textbf{Proof:} Let $\left[ \begin{array}{c}
\vec y\\
\vec z \end{array} \right]\in \mathrm H.$ Then $\left<B \left[ \begin{array}{c}
\vec y\\
\vec z \end{array} \right], \left[ \begin{array}{c}
\vec y\\
\vec z \end{array} \right] \right>_{\mc H}=  k_2 |z_3(L)|^2.$
The compactness follows from the fact that $M^{-1}$ is a canonical isomorphism from $\mathrm H$ to $\mathrm H',$ and the fact that $B$ is a rank-one operator, hence compact from $\mathrm H$ to $\mathrm H'.$ $\square$

\subsection{Description of ${\rm Dom}(\mc A)$}

\begin{prop}\label{prop-dom}Let  $\vec u=(\vec y, \vec z)^{\rm T}\in \mc H.$ Then $\vec u \in {\rm Dom}(\mc A)$  if and only if the following conditions hold:
\begin{eqnarray}
\nonumber && \vec y\in (H^2(0,L)\cap H^1_L(0,L))^2 \times (H^3(0,L) \cap H^2_L(0,L)),\quad \vec z\in V ~{\rm such ~ that}~(y_1)_{x}=(y_2)_{x}=(y_4)_{xx}\left. \right|_{x=L}=0.
\end{eqnarray}
Moreover, the resolvent of $\mc A$ is compact in the energy space $\mc H.$
\end{prop}
\vspace{0.1in}

\textbf{Proof:} Let $\vec {\tilde u}= \left( \begin{array}{c}
\vec {\tilde y} \\
\vec {\tilde z}\\
 \end{array} \right)\in \mc H$ and $\vec {u}= \left( \begin{array}{c}
\vec {y} \\
\vec {z}\\
 \end{array} \right)\in {\rm Dom}(\mc A)$ such that $\mc A \vec u=\vec {\tilde u}.$ Then we have
 $$ -\vec z= \vec {\tilde y} \in V, \quad   A\vec y + D_0 z=M\vec{\tilde z},$$ and therefore,
\begin{eqnarray}\label{sal}\left< \vec y, \vec \varphi\right>_{V}= \left< \vec {\tilde z}, \vec \varphi \right>_{H} ~~{\rm for ~~all ~~} \vec \varphi \in V.
\end{eqnarray}
Let  $\vec \psi=[\psi_1,\psi_2,\psi_3]^{\rm T} \in (C_0^\infty(0,L))^4.$  We define $\varphi_i=\psi_i$ for $i=1,2,$ and $\varphi_3=\int_0^x \psi_3 (s)ds.$ Since $\vec \varphi\in V,$ inserting $\vec \varphi$ into the above equation yields
\begin{eqnarray} \nonumber &&\int_0^L \left\{  -\alpha^1 h_1  (y_1)_{xx} \bar \psi_1  - \left(\alpha^3 h_3  (y_2)_{xx}+\frac{\gamma^2 h_3}{\ep3} (P_\xi (y_2)_x)_x\right) \bar \psi_2- K_2 (y_3)_{xxx} \bar \psi_3   \right\}~dx\\
\nonumber &&  + \frac{G_2}{h_2}  (-y_1+y_2 + H(y_3)_x)(-\bar \psi_1+\bar \psi_2 + H(\bar \psi_3)_x)+s_1 (z_1)_x(L)(\psi_1)_x(L)  + \frac{\gamma s_3}{\ep3} z_2(L)\psi_2(L)+ k_1(z_3)_{x}(L) \psi_3(L) \\
  \nonumber &&= \int_0^L \left\{ \left(\int_1^x m \tilde z_4  ds + K_1 (\tilde z_4)_x\right) \bar \psi_4 +\rho_1 h_1 \tilde z_1 \bar \psi_1+ \rho_3 h_3 \tilde z_2 \bar \psi_2 + \mu h_3 \tilde z_3 \bar \psi_3 \right\}~dx\\
 \end{eqnarray}
for all $\vec \psi\in (C_0^\infty(0,L))^3.$ Therefore it follows that $\vec y\in (H^2(0,L)\cap H^1_L(0,L))^2 \times (H^3(0,L) \cap H^2_L(0,L)).$

Next let $\vec \psi \in \mathrm H.$  We define
\begin{eqnarray} \label{dumber}\varphi_i=\int_0^x \psi_i(s)ds,\quad i=1,\ldots,3.
\end{eqnarray}
Obviously $\vec \varphi \in \rm V. $ Then plugging (\ref{dumber}) into (\ref{sal}) yields
\begin{eqnarray}&& \nonumber  0=(\alpha^1 h_1  (y_1)_{x}(L) + s_1 z_1(L)) \bar \psi_1(L)+\alpha^3 h_3  (y_2)_{x}(L)+ \frac{\gamma^2 h_3}{\ep3} P_\xi (y_2)_x(L) \\
\nonumber &&+ \frac{\gamma s_3}{\ep3} z_2(L) \bar \psi_2(L)+ (k_1 (y_3)_{xx}(L) + k_1 (z_3)_x(L)) (\bar \psi_3)_x(L)
\end{eqnarray}
for all $\psi\in \mathrm H.$ Hence,
\begin{eqnarray}\nonumber &&\alpha^1 h_1  (y_1)_{x}(L) + s_1 z_1(L)=\alpha^3 h_3  (y_2)_{x}(L) + \frac{\gamma^2 h_3}{\ep3} P_\xi (y_2)_x(L)+ \frac{\gamma s_3}{\ep3} z_2(L)= k_1 (y_3)_{xx}(L) + k_1 (z_3)_x(L)=0.
\end{eqnarray}

Now let  $\vec y=\left[ \begin{array}{c}
\vec y_1 \\
\vec y_2 \end{array} \right]\in \mc {\rm Dom} (\mc A )$  and $\vec z=\left[ \begin{array}{c}
\vec z_1 \\
\vec z_2 \end{array} \right] $such that $(I+\mc A) \vec y=\vec z.$ By Proposition \ref{prop-dom} and  Lemma \ref{skew},  the compactness of the resolvent follows. $\square$

\begin{lem} \label{xyz} The eigenvalue problem
\begin{eqnarray}
 \label{dbas-eig} &&\left\{
  \begin{array}{ll}
   \alpha^1 h_1 z^1_{xx} - G_2 \phi^2 = \lambda^2 \rho_1 h_1 z^1,& \\
 \alpha^3 h_3 z^3_{xx} + \frac{\gamma^2 h_3}{\ep3} (P_\xi(z_3)_x)_x + G_2  \phi^2 = \lambda^2 \rho_3 h_3 z^3   & \\
    - K_2 u_{xxxx} + G_2 H \phi^2_x=\lambda^2 (m u  - K_1 u_{xx}),&
 \end{array} \right.
\end{eqnarray}
with the overdetermined boundary conditions
\begin{eqnarray}
 \nonumber &&  u(0)=u_x(0)=z^1(0)=z^3(0)=z^1(L)=z^1_x (L)=z^3(L)=z^3_x(L)=0,\\
  \label{d-son-eig} && u(L)=u_x(L)=u_{xx}(L)=u_{xxx}(L)=0
\end{eqnarray}
has only the trivial solution.
\end{lem}
\vspace{0.1in}

\textbf{Proof:} Now multiply the equations in (\ref{dbas-eig}) by $-x\bar u_x+3\bar u,$ $ -x\bar z^3_x +2\bar z^3,$ and $-x \bar z^1_x+2\bar z^1,$  respectively, integrate by parts on $(0,L),$

\begin{eqnarray}
\nonumber && 0=\int_0^L\left\{ -\alpha^1 h_1|z^1_x|^2 +(x\bar z^3_{xx} -\bar z^3_x) \left(\alpha^3 h_3 z^3_x +\frac{\gamma^2 h_3}{\ep3} P_\xi z^3_x\right)  -3\rho_1 h_1 \lambda^2|z^1|^2 -3\rho_3 h_3  \lambda^2 |z^3|^2 -4m \lambda^2 |u|^2  \right. \\
  \nonumber && -  2 K_1 \lambda^2|u_x|^2  - G_2 h_2  \phi^2 (x \bar \phi^2_x)-3G_2 h_2|\phi^2|^2-K_2 \bar u_{xxxx} (x u_x) + \alpha^1_1 h_1 \bar z^1_{xx} (x z^1_x)  -\rho_1h_1 \lambda^2\bar z^1 (xz^1_x)\\
\label{CH3-mult1-20} && \left. -\rho_3h_3 \lambda^2 \bar z^3 (x z^3_x)  -\lambda^2 (m\bar u-K_1 u_{xx})(x u_x)\right)~dx\quad\quad
 \end{eqnarray}
Now consider the conjugate eigenvalue problem corresponding to (\ref{dbas-eig})-(\ref{d-son-eig}). Now multiply the equations in the conjugate problem by $-x u_x-2 u,$ $-x z^1_x-3 z^1,$ and $ -x z^3_x-3 z^3,$ respectively, integrate by parts on $(0,L),$ and add them up:
\begin{eqnarray}
\nonumber && 0=\int_0^L\left\{ 3\alpha^1 h_1 |z^1_x|^2   -x z^3_{x} \left(\alpha^3 h_3 \bar z^3_{xx} +\frac{\gamma^2 h_3}{\ep3} (P_\xi \bar z^3_x)_x\right) +3 z^3_{x} \left(\alpha^3 h_3 \bar z^3_{x} +\frac{\gamma^2 h_3}{\ep3} P_\xi \bar z^3_x\right) \right. \\
  \nonumber && +3\bar\lambda^2\rho_1 h_1 |z^1|^2 +3\rho_3 h_3 \bar \lambda^2|z^3|^2  +2 m \bar \lambda^2 |u|^2 +  2 K_1 \bar \lambda^2|u_x|^2 +2K_2 |u_{xx}|^2  + G_2 h_2 \bar \phi^2 (x \phi^2_x)+3G_2 h_2|\phi^2|^2\\
\label{CH3-mult1-21} && \left.+K_2 \bar u_{xxxx} (x u_x) - \alpha^1_1 h_1 \bar z^1_{xx} (x z^1_x)   +\rho_1h_1 \bar \lambda^2\bar z^1 (xz^1_x)+\rho_3h_3 \bar \lambda^2 \bar z^3 (x z^3_x)  +\bar \lambda^2 (m\bar u-K_1 u_{xx})(x u_x)\right)~dx\quad\quad
 \end{eqnarray}

Adding (\ref{CH3-mult1-20}) and (\ref{CH3-mult1-21}) yields,
\begin{eqnarray}
\nonumber && 0=\int_0^L\left\{ 2\alpha^1 h_1 |z^1_x|^2  +2\alpha^3 h_3 |z^3_x|^2 + \frac{2\gamma^2 h_3}{\ep3} (P_\xi z^3_x\cdot \bar z^3_x) +\frac{\gamma^2 h_3}{\ep3} x\left( -z^3_{x} (P_\xi \bar z^3_x)_x +\bar  z^3_{xx} P_\xi \bar z^3_x\right) \right. \\
  \nonumber && +3\bar\rho_1 h_1 (-\lambda^2 + \bar \lambda^2) |z^1|^2 +3\rho_3 h_3 \bar  (-\lambda^2 + \bar \lambda^2)|z^3|^2  + (-4m \lambda^2 + 2 m \bar \lambda^2) |u|^2 +  2K_1(- \lambda^2 + \bar \lambda^2)|u_x|^2 \\
\label{CH3-mult1-22} && \left.+2K_2 |u_{xx}|^2  + G_2 h_2 \left(- \phi^2 (x \bar \phi^2_x) + \bar \phi^2 (x \phi^2_x)\right)+ (-\lambda^2+\bar \lambda^2) (m\bar u-K_1 u_{xx})(x u_x)\right)~dx\quad\quad
 \end{eqnarray}

Finally, adding (\ref{CH3-mult1-20}) and (\ref{CH3-mult1-21}),considering only the real part of the expression above and  all eigenvalues are
located on the imaginary axis, i.e. $\lambda=\mp \imath \nu,$ yields
\begin{eqnarray}
\nonumber && \int_0^L\left( K_2 |u_{xx}|^2+ m \nu^2|u|^2+ \alpha^1 h_1 |z^1_x|^2  +\left(\alpha^3 h_3 z^3_x + \frac{\gamma^2 h_3}{\ep3} (P_\xi z^3_x)\right)\cdot \bar z^3_x \right.\\
\label{damn1}&& \quad \left.+\frac{\gamma^2 h_3}{2\ep3} x\left( -z^3_{x} (P_\xi \bar z^3_x)_x +\bar  z^3_{xx} P_\xi \bar z^3_x\right) \right)~dx=0
 \end{eqnarray}

Now we should get rid of the last term which survived through the calculations. Let $P_\xi z^3_x=\eta.$ Then $\eta-\xi \eta_{xx}=z^3_x,$ and
\begin{eqnarray}
\nonumber  &&\frac{\gamma^2 h_3}{2\ep3} \int_0^L x\left[ -z^3_{x} (P_\xi \bar z^3_x)_x +\bar  z^3_{xx} P_\xi \bar z^3_x\right]~dx = \frac{\gamma^2 h_3}{2\ep3} \int_0^L x\left[ (xz^3_{xx}+z^3_x)(P_\xi \bar z^3_x) +\bar  z^3_{xx} P_\xi \bar z^3_x\right]~dx\\
\nonumber && \quad = \frac{\gamma^2 h_3}{2\ep3} \int_0^L \left[ (x(\eta_x-\xi \eta_{xxx})+(\eta-\xi\eta_{xx}))\bar \eta +x(\bar\eta_x-\xi\bar \eta_{xxx}) \eta\right]~dx\\
\nonumber && \quad =  \frac{\gamma^2 h_3}{2\ep3} \int_0^L \left[ x (\eta\cdot \bar\eta)_x + |\eta|^2 + \xi |\eta_x|^2 + \xi \left(\eta_{xx} (\bar\eta+x\bar\eta_x) + \bar\eta_{xx}(\eta+\eta_{xx})\right)\right]~dx\\
\label{damn2} && \quad = -\frac{\gamma^2 h_3\xi}{2\ep3}\int_0^L |\eta_x|^2 dx,
\end{eqnarray}
and
\begin{eqnarray}
\label{damn3} \frac{\gamma^2 h_3}{\ep3} \int_0^L  (P_\xi z^3_x) \cdot \bar z^3_x ~dx&=&  \frac{\gamma^2 h_3}{\ep3} \int_0^L  \eta (\bar \eta -\xi \bar \eta_{xx}) ~dx=\frac{\gamma^2 h_3}{\ep3} \int_0^L  (|\eta|^2 + \xi |\eta_x|^2) ~dx.
 \end{eqnarray}
Plugging (\ref{damn2}) and (\ref{damn3}) in (\ref{damn1}) gives
\begin{eqnarray}
\nonumber && \int_0^L\left( K_2 |u_{xx}|^2+ m \nu^2|u|^2+ \alpha^1 h_1 |z^1_x|^2  +\left(\alpha^3 h_3 z^3_x + \frac{\gamma^2 h_3}{\ep3} (P_\xi z^3_x)\right)\cdot \bar z^3_x \right.\\
\label{damn4} && \quad \left.+ \frac{\gamma^2 h_3}{\ep3}  (|P_\xi z^3_x|^2 + \frac{\xi}{2} |(P_\xi z^3_x)_x|^2) \right)~dx=0.
 \end{eqnarray}
This together with (\ref{d-son-eig}) imply that $u=z^1=z^3\equiv 0$ by  (\ref{d-son-eig}). This also covers the case $\lambda\equiv 0. \square$

Now we  consider the decomposition $\mc{A}+B=(\mc{A}_d+B)+ \mc {A}_{\phi}$ of the semigroup generator of the original problem (\ref{defA}) where $\mc{A}_d+B$ is the semigroup generator of the decoupled system, i.e. $\phi^2\equiv 0$ in (\ref{d4-non})-(\ref{d-son-non}),
\begin{eqnarray}
 \label{dbas-dcpld} &&\left\{
  \begin{array}{ll}
   m \ddot w - K_1 \ddot{w}_{xx}+K_2w_{xxxx} =0,&\\
    \rho_1h_1 \ddot v^1  -\alpha^1 h_1 v^1_{xx}  =0,   &  \\
  \rho_3h_3 \ddot v^3       -\alpha^3 h_3 v^3_{xx}    = 0,   &
 \end{array} \right.
\end{eqnarray}
with the boundary and initial conditions
 \begin{eqnarray}
&&  \begin{array}{ll}
\nonumber   \left\{v^1, v^3, w, w_x  \right\}_{x=0}=0, ~  \left\{\alpha^1 h_1v^1_x\right\}_{x=L}=-s_1 \dot v^1(L),~~  \left\{ \alpha^3 h_3 v^3_x  \right\}_{x=L}=-s_3 \dot v^3(L,t),&\\
  \end{array} \\
         &&  \begin{array}{ll}
\nonumber  \left\{K_2 w_{xx}  \right\}_{x=L} =-k_1\dot w_x(L,t) , ~  \left\{K_1 \ddot w_x -K_2 w_{xxx}  \right\}_{x=L} =k_2 \dot w(L,t),& t\in \mathbb{R}^+ \\
  \end{array} \\
&&  \begin{array}{ll}
\label{d-son-dcpld}  (v^1,v^3, w, \dot v^1, \dot v^3, \dot w)(x,0) =(v^1_0, v^3_0, w_0, v^1_1, v^3_1, w_1),& x \in [0,L].
  \end{array}
\end{eqnarray}
The operator $\mc{A}_{\phi}:\mc{H}\to\mc{H}$ is the coupling between the layers defined as the following
\begin{eqnarray}\label{opB}\mc {A}_{\phi}{\bf y}=
 \left( \begin{array}{c}
 0_{3\times 1} \\
 \mc  M^{-1} \left(H G_2 ~\phi^2_x\right)   \\
\frac{G_2}{h_1 \rho_1}\phi^2\\
\frac{{\gamma^2 h_3}}{{\ep3}}~(P_{\xi} v^3_x)_x-\frac{G_2}{h_3 \rho_3}\phi^2\\
   \end{array} \right).
 \end{eqnarray}
 where ${\bf y}=(w,u^1, u^3, \tilde w, \tilde v^1, \tilde v^3)$ and  $\phi^2=\frac{1}{h_2}\left(-u^1+u^3 + H u_x\right).$
 Let $E_d(t)$ be natural energy corresponding to the system (\ref{dbas-dcpld})-(\ref{d-son-dcpld}), i.e. (\ref{Energy-nat-non}) without the $P_\xi$ and $\phi^2$ terms.

 \begin{thm} Let ${\mc A}_d+B$ be the infinitesimal generator of the semigroup corresponding to the solutions of (\ref{dbas-dcpld})-(\ref{d-son-dcpld}). Then the semigroup $\{e^{({\mc A}_d+B) t}\}_{t\ge 0}$ is exponentially
stable in $\mc H$.
 \end{thm}

 \textbf{Proof:}  Note that the equations in (\ref{dbas-dcpld}) are completely decoupled. The exponential stability of the semigroup $e^{({\mc A}_d+B)t}$ follows from the exponential stability of wave equations  \cite{O-Hansen4} and the Rayleigh beam equation \cite{B-Rao}.

\begin{lem} \label{compact} The operator $A_{\phi}:\mc{H}\to \mc{H}$ defined in (\ref{opB}) is compact.
\end{lem}

\vspace{0.1in}
{\bf Proof:}  When $(w,v^1, v^3, \tilde w, \tilde u^1, \tilde u^3)\in \mathrm H,$ we have $w\in H^2_L(0,L)$ and $v^1, v^3 \in (H^1_L(0,L))^{2},$ and therefore $\phi^2 \in H^1_L(0,L).$ Since $\mc M: H^2_L(0,L)\to \mathrm L^2(0,L)$ remains  an isomorphism, the terms in (\ref{opB}) satisfy
\begin{eqnarray}
\label{dumb5} && \mc M^{-1} \left(\phi^2_x\right)  \in H^2_L(0,L), ~~\phi^2 \in H^1_L(0,L),~~ (P_\xi v^3_x)_x \in H^1_0(0,L)
 \end{eqnarray}
and  $H^2_L(0,L)\times  H^1_L(0,L)$ is compactly embedded  in $H^1_L(0,L)\times (\mathrm L^2(0,L))^{2},$ and $H^1_0(0,L)$ is compactly embedded in $L^2(0,L).$ Hence the operator $A_{\phi}$  is compact in $\mc H$.

 \begin{thm} \label{main-thm} Then the semigroup $\{e^{({\mc A}+B) t}\}_{t\ge 0}$ is exponentially
stable in $\mc H.$
 \end{thm}

 \textbf{Proof:} The semigroup $\mc A+B={\mc A}_d +B + {\mc A}_{\phi}$ is strongly stable on $\mc H$ by Lemma \ref{xyz}, and the operator ${\mc A}_{\phi}$ is a compact in $\mc H$ by
Lemma \ref{compact}. Therefore, since the semigroup generated by $({\mc A}_d +B + {\mc A}_{\phi})-{\mc A}_{\phi}$ is uniformly
exponentially stable in $\mc H,$  the semigroup $\mc A=({\mc A}_d + B+ {\mc A}_{\phi})$ is uniformly
exponentially stable in ${\mc H}$ by e.g., the perturbation theorem of \cite{Trigg}.

The main difference between the voltage and charge controlled models is the existence of the $P_\xi$ term in the charge controlled model. By using the same argument above (by taking $P_\xi\equiv 0$), the voltage controlled model can be easily shown to be exponentially stable with the same feedback controller (\ref{feedback}).

\begin{thm}  Then the semigroup $\{e^{({\mc A}+B) t}\}_{t\ge 0}$ corresponding to (\ref{d4-non}) ($P_\xi\equiv 0$) with (\ref{feedback}) is exponentially
stable in $\mc H.$
 \end{thm}

 Several remarks are in order:

\begin{rmk} (i) The number of controllers for the bending motion may be reduced to one by taking $k_2=0.$ However, it is worthwhile to note that the multiplier used in \cite{Rao} does not work in the proof of Lemma \ref{xyz}. For that reason, the classical four boundary conditions for $u$ at $x=L$ in (\ref{d-son-eig}) are still needed due to the existence of the coupling $\phi^2$ in (\ref{dbas-eig}). Removing $u(L)=0,$ and proving the same result in Lemma \ref{xyz} is an open question.

\noindent (ii) If the moment of inertia term in (\ref{d4-non})-(\ref{d-son-non}) is set to zero, $K_1\equiv 0,$ the bending motion is described by the Euler-Bernoulli-type equation. The exponential stabilizability (\ref{d4-non})-(\ref{d-son-non}) with $s_2\equiv 0,$ i.e. one controller for each equation, is still an open problem.
\end{rmk}
\section{Generalization to the Multilayer ACL beam case}
\begin{figure}
\begin{center}
\includegraphics[height=5.5cm]{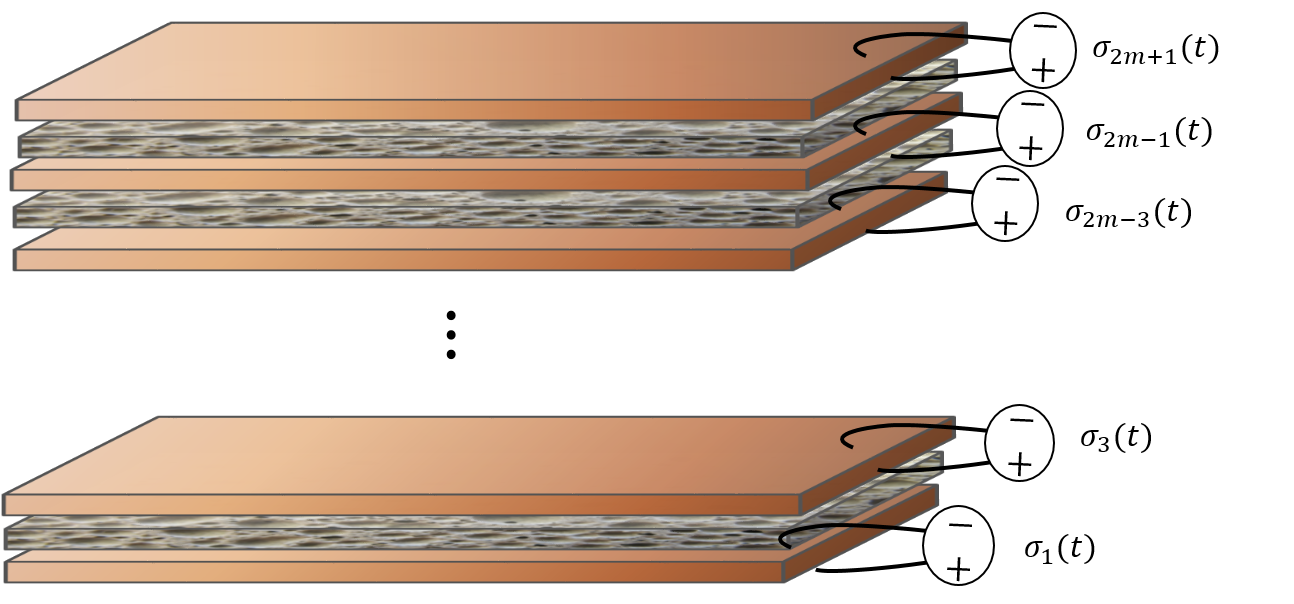}    
\caption{The composite consists of joint ACL layers. Odd layers are all piezoelectric and are actuated by a charge source.}  
\label{ACL-mult}                                 
\end{center}                                 
\end{figure}

By using the same methodology in Section \ref{modeling}, we can obtain the model for the multilayer generalization of the ACL sandwich beams as in Figure \ref{ACL-mult}.
The equations of motion are found to be
\begin{equation}\left\{ \begin{array}{l}
   {\bf{h}}_{\mc O} {\bf{p}}_{\mc O} {\ddot y}_{\mc O} -{\bf{h}}_{\mc O} {\bf{\alpha}}_{\mc O}  ({y}_{\mc O})_{xx} -{\bf{\gamma}}_{\mc O}^2{\bf{\varepsilon}}_{\mc O}^{-1} {\bf{h}}_{\mc O} ~(P_{\xi_{\mc O}} ({y}_{\mc O})_{x})_x+ {\bf{B}}^{\rm T}  {\bf{ G}}_E \psi_E  = \delta(x-L){\bf \gamma}_{\mc O}{\bf \varepsilon}_{\mc O}^{-1}\sigma_{\mc O}(t), \\
m\ddot w -K_1 \ddot w_{xx} +  K_2   w_{xxxx} -   N^{\rm T} {\bf{h}}_E {\bf{ G}}_E \psi_E = 0 ~~~ {\rm{in}}~~ (0, L)\times \mathbb{R}^+ \\
 {\text{where}}~~  {\bf{B}} { y}_{\mc O}={\bf{h}}_E \phi_E-{\bf{h}}_E  N w',
 \end{array} \right.
\label{maincont}
\end{equation}
with clamped-free boundary conditions  and initial conditions
\begin{eqnarray}\label{bdryback10}
 &  \left\{ \begin{array}{l}
  \left\{{y}_{\mc O}, w, w_x  \right\}_{x=0}=0 ~~ {\rm{on}}~ \mathbb{R}^+,\label{bdrycont3}\\
   \left\{\alpha_{\mc O} h_{\mc O}y_{\mc O} + {\bf{\gamma}}_{\mc O}^2{\bf{\varepsilon}}_{\mc O}^{-1} {\bf{h}}_{\mc O} ~P_{\xi_{\mc O}} ({y}_{\mc O})_{x}\right\}_{x=L}=0,~~~ {\rm{on}}~~ \mathbb{R}^+,\label{bdrycont4}\\
   \left\{K_2 w_{xx}  \right\}_{x=L} =-M(t) , ~  \left\{K_1 \ddot w_x -K_2 w_{xxx}  \right\}_{x=L} = g(t), ~~ {\rm{on}}~ \mathbb{R}^+ ,\label{bdrycont5}\\
w(x,0)=w^0(x),  ~\dot w(x,0)=w^1(x), ~{y}_{\mc O}(x,0)= {y}^0_{\mc O}, ~ {\dot y}_{\mc O}(x,0)= {y}^1_{\mc O} ~ {\rm{on}}~~ (0,L).  \label{initialcont}
   \end{array} \right. &
\end{eqnarray}

The model (\ref {maincont})  consists of  $2m+1$ alternating stiff and complaint (core) layers, with piezoelectric layers on outside.   The piezoelectric layers have odd indices $1,3,\ldots 2m+1$ and the even layers have even indices $2,4,\ldots 2m$.

 In the above, $m,K_1, K_2$ are \emph{positive} physical constants, $w$ represents the transverse displacement, $\phi^i$ denotes the shear angle in the $i^{\rm{th}}$ layer, $\phi_E=[\phi^2,\phi^4,\ldots,\phi^{2m} ]^{\rm T},$ $y^i$ denote the longitudinal displacement
along the center of the   $i^{\rm{th}}$   layer, and $y_{\mc O}=[y^1,y^3, \ldots, y^{2m+1}]^{\rm T}.$  Define the following $n\times n$ diagonal matrices
\begin{eqnarray}
\nonumber &{\bf{p}}_{\mc O}={\rm{diag}}~ (\rho_1,  \ldots,\rho_{2m+1}),  ~~{\bf{h}}_{\mc O}={\rm{diag}}~(h_1, \ldots, h_{2m+1}), ~~{\bf{h}}_E={\rm{diag}}~(h_2,\ldots, h_{2m}),&\\
\nonumber & {\bf{\varepsilon}}_{\mc O}={\rm{diag}}~(\varepsilon_3^1, \ldots, \varepsilon_3^{2m+1}), ~~{\bf{\xi}}_O={\rm{diag}}~(\xi_1,\ldots, \xi_{2m+1})~~{\bf{\alpha}}_{\mc O}={\rm{diag}}~(\alpha_1, \ldots, \alpha_{2m+1}),&\\
\nonumber & {\bf{ G}}_E={\rm{diag}}~( G_2, \ldots, G_{2m}), ~~{\bf{\sigma}}_{\mc O}={\rm{diag}}~(\sigma_1, \ldots, \sigma_{2m+1})&
\end{eqnarray}
where $h_i, \rho_i, E_i,$ are positive and denote the thickness, density, and Young's modulus, respectively.
Also $G_i\ge 0$ denotes shear modulus of the $i^{\rm{th}}$ layer, and $ {\tilde{ G}}_i\ge 0$ denotes coefficient for damping in the corresponding compliant layer.

The vector $ N$ is defined as
$N={\bf{h}}_E^{-1}{\bf{A }}{\bf{h}}_{\mc O}  \vec 1_{\mc O} +  \vec 1_E$
where ${\bf{A}}=(a_{ij})$  and ${\bf{B}}=(b_{ij})$ are the $m\times(m+1)$ matrices
$$a_{ij}  = \left\{ \begin{array}{l}
1/2,~~{\rm{  if   }}~~j = i~~{\rm{  or }}~~j = i + 1 \\
~~0,\quad{\rm{   otherwise}} \\
\end{array} \right., ~~b_{ij}=\left\{ \begin{array}{l}
(-1)^{i+j+1},~~{\rm{  if   }}~~j = i~~{\rm{  or }}~~j = i + 1 \\
~~0, \quad\quad\quad\quad {\rm{    otherwise}} \\
\end{array} \right. $$
and $ \vec 1_{\mc O}$ and $\vec 1_E$ denote the vectors with all entries $1$ in $\mathbb{R}^{m+1}$ and $\mathbb{R}^{m},$ respectively.

In the above, the operator $P_{\xi_{\mc O}}$ is defined by
\begin{eqnarray}\label{Lgamma-mat}&&P_{\xi_{\mc O}}:={\rm{diag}}~(P_{\xi_1}, \ldots, P_{\xi_{2m+1}})\end{eqnarray}
where $$P_{\xi_i}=\left(-\xi_i D_x^2+I\right)^{-1}, \quad \xi_i:=\frac{\e1^i h_i^2}{12\ep3^i}, ~~~i=1,3,\ldots, 2m+1.$$

Let ${\bf{s}}_O={\rm{diag}}~(s_1,\ldots, s_{2m+1}),~$ $s_i,k_1, k_2>0,$  and
 \begin{eqnarray}\label{feedback-mult} {\bf \sigma}_{\mc O}^i=-{\bf s}_{\mc O}\dot {y}_{\mc O}(L,t),~~M(t)= k_1\dot w_x(L,t),~~ g(t)=k_2 \dot w(L).
\end{eqnarray}
The well-posedness of the model (\ref{maincont})-(\ref{initialcont}) with (\ref{feedback-mult}) can be established in a similar fashion as in Section \ref{Section-stab}. The following result follows:
\begin{thm} The system (\ref{maincont})-(\ref{initialcont}) with (\ref{feedback-mult}) is exponentially stable.
\end{thm}\\
\textbf{Proof:} The proof is analogous to the proof of Theorem \ref{main-thm}.

\begin{rmk} The equations of motion for a voltage controlled multilayer ACL beam can be obtained similarly. One can also consider a hybrid type \emph{hybrid} multilayer ACL beam model for which some layers are actuated by voltage sources, and the rest are actuated by charge sources. The equations of motion can be also obtained in a similar fashion.
\end{rmk}
 \section{Conclusion and Future Research}
In this paper, it is shown that in the case of electrostatic assumption, charge, current, or voltage actuated ACL beams can be uniformly exponentially stabilized by using mechanical feedback controllers. This is not the case once we consider the quasistatic or fully dynamic approaches where the magnetic effects are accounted for. For example, the voltage-actuated ACL beam model obtained by the fully dynamic or quasistatic approaches is shown to be not uniformly exponentially stabilizable \cite{Ozkan3}. The polynomial stability result obtained in \cite{Ozkan1} for certain combinations of material parameters on a more regular space than the natural energy space is applicable to this problem yet it is still an open problem. Moreover, the fully dynamic model  for current or charge actuation is obtained in \cite{Ozkan4}, yet the stabilization problem  is currently open and under investigation.

Even though the quadratic-through-thickness assumption for the electric potential in the case of charge or current actuation makes the beam stiffer than voltage actuated ACL beam, this effect is not observed for higher order eigenvalues due to the compactness of the operator $P_\xi$. 

%



\bibliographystyle{spiebib} 

\bibliographystyle{plain}        

\begin{thebibliography}{99}     


\bibitem{Baz} A. Baz, {``Boundary Control of Beams Using Active Constrained Layer Damping,"} {\sl J. Vib. Acoust.} {\bf 119(2)}, 166--172 (1997).

\bibitem{Cao-Chen} Y. Cao, X.B. Chen, {``A Survey of Modeling and Control Issues for Piezo-electric Actuators,"} {\sl Journal of Dynamic Systems, Measurement, and Control } {\bf 137(1)}, 014001 (2014).



\bibitem{Review} S Devasia, E. Eleftheriou, S. O. Reza Moheimani, {``A Survey of Control Issues in Nanopositioning,"} {\sl IEEE Transactions on Control Systems Technology} {\bf 15(5)}, 802--823 (2007).

\bibitem{Ditaranto} R.A. DiTaranto, { ``Theory of vibratory bending for elastic and viscoelastic layered finitelength
beams,"} {\sl J. Appl. Mech.} {\bf 32}, 881--886 (1965).



\bibitem{F} R.H. Fabiano, S.W. Hansen, { ``Modeling and analysis of a three--layer
damped sandwich beam,"} {\sl Dynamical systems and differential
equations,  Discrete Contin. Dynam. Systems} \textbf{Added Volume}, 143--155 (2001).

\bibitem{Hansen} S.W.  Hansen, {``Analysis of a Plate with a Localized Piezoelectric Patch,"}   {\sl The Proceedings of the IEEE Conference on Decision \& Control,} Tampa, Florida, 2952-2957 (1998).

\bibitem{Hansen3}{  S.W. Hansen,}
{ ``Several Related Models for Multilayer Sandwich Plates,"}
{\sl  Mathematical Models \& Methods in Applied
Sciences}  {\bf 14(8)}, 1103-1132 (2004).

 \bibitem{Lee} P.C.Y. Lee, {`` A variational principle for the equations of piezoelectromagnetism in elastic dielectric crystals,"} Journal of Applied Physics {\bf{(69(11)}}, 7470--7473 (1991).

\bibitem {Main1} J.A. Main and E. Garcia, {``Design impact of piezoelectric actuator
nonlinearities,"} {\sl Journal of Guidance, Control, and Dynamics} {\bf 20(2)}, 327-–332 (1997).

\bibitem {Main2} J.A. Main, E. Garcia and D.V. Newton, {``Precision position control of
piezoelectric actuators using charge feedback,"} {\sl  Journal of Guidance, Control, and Dynamics}  {\bf 18(5)},  1068–-1073 (1995).

 \bibitem{Mead} D.J. Mead and S. Markus, {``The forced vibration of a three-layer, damped sandwich beam
with arbitrary boundary conditions,"} {\sl J. Sound Vibr.} {\bf 10}, 163--175 (1969).


\bibitem{M-F} S.O.R. Moheimani, A.J. Fleming, {[Piezoelectric transducers for vibration control and damping],} {Springer-Verlag} (2006).
\bibitem{O-M1}  K.A. Morris, A.\"{O}. \"{O}zer, {{``Modeling and stabilizability of voltage-actuated piezoelectric beams with magnetic effects,"}}  {\sl SIAM J. Control Optim.} {\bf 52(4)}, 2371--2398 (2014).

\bibitem{O-M3}  K.A. Morris, A.\"{O}. \"{O}zer, {``Comparison  of stabilization of  current-actuated  and voltage-actuated piezoelectric beams,"}{\sl The $53^{\rm  rd}$ Proceedings of the IEEE Conf. on Decision \& Control,} Los Angeles, California, USA, 571--576 (2014).



\bibitem{Ozkan1} A.\"{O}. \"{O}zer, {\newblock{``Further stabilization and exact observability results for voltage-actuated piezoelectric beams with magnetic effects,"}} {\sl Mathematics of Control, Signals, and Systems} {\bf 27(2)}, 219--244 (2015).

\bibitem{Ozkan2} A.\"{O}. \"{O}zer, {\newblock{``Modeling and well-posedness results for active constrained layered (ACL) beams with/without magnetic effects,"}} accepted by  The Proceedings of the American Control Conference, (2016).
\bibitem{Ozkan3} A.\"{O}. \"{O}zer, {\newblock{``Modeling and control results for an active constrained layered (ACL) beam  actuated by two voltage sources with/without magnetic effects,"}} submitted, arXiv:1511.05907.
    \bibitem{Ozkan4} A.\"{O}. \"{O}zer, {\newblock{``    Modeling and semigroup formulation of charge or current-controlled active constrained layer (ACL) beams; electrostatic, quasi-static, and fully-dynamic cases,"}} submitted.
\bibitem{O-Hansen3} {A.\"{O}. \"{O}zer, and S.W. Hansen},  \newblock{``Exact boundary controllability results for a multilayer Rao-Nakra sandwich beam,"} {\sl SIAM J. Control Optim.} {\bf 52(2)}, 1314--1337 (2014).

\bibitem{O-Hansen4} {A.\"{O}. \"{O}zer, and S.W. Hansen}, \newblock{``Uniform stabilization of a multi-layer Rao-Nakra sandwich beam,"}  {\sl Evolution Equations and Control Theory} {\bf 2(4)}, 195--210 (2013).

 \bibitem{B-Rao} B. Rao, {``A compact perturbation method for the boundary stabilization of the Ragleigh beam equation,"}  {\sl Appl. Math. Optim.} {\bf 33(3)},  253--264 (1996).
\bibitem {Rao} Y.V.K.S. Rao and B.C. Nakra, {``Vibrations of unsymmetrical sandwich beams and plates with viscoelastic cores,"} {\sl J. Sound Vibr.} {\bf 34(3)}, 309--326 (1974).

\bibitem{Rogacheva} N. Rogacheva, {[The Theory of Piezoelectric Shells and Plates]}, Boca Raton, FL: CRC Press (1994).

\bibitem{Smith} R.C. Smith, { [Smart Material Systems: Model Development]}, Society for Industrial and Applied Mathematics (2005).

\bibitem{Stanway} R. Stanway, J.A. Rongong, N.D. Sims, {``Active constrained-layer damping: a state-of-the-art review,"}  {\sl Automation \& Control Systems }{\bf 217(6)}, 437--456 (2003).
\bibitem {Sun} C.T. Sun and Y.P. Lu, {[Vibration Damping of Structural Elements],} Prentice Hall (1995).
\bibitem{Tabesh} A. Tabesh and L.G. Fr\'{e}chette, {``An improved small-deflection electromechanical model for piezoelectric bending beam actuators and energy harvesters,"} {\sl J. Micromech. Microeng.} {\bf 18-104009
}, 1-12 (2008)
\bibitem{Trigg}  R. Triggiani, {``Lack of uniform stabilization for noncontractive semigroups under compact perturbation,"} {\sl Proc. Amer. Math. Soc.}  {(\bf 105)}, 375-383 (1989).

\bibitem{Trindade} M. Trindade and A. Benjendou, { ``Hybrid Active-Passive Damping Treatments Using Viscoelastic and Piezoelectric Materials:Review and Assessment,"} {\sl Journal of Vibration and Control} {\bf 8(6)}, 699--745 (2002).


\bibitem{Yan} M.-J. Yan and E.H. Dowell, {``Governing equations for vibrating constrained-layer damping
sandwich plates and beams,"} {\sl J. Appl. Mech.} {\bf 39}, 1041--1046 (1972).
\end{thebibliography}

\end{document}